\documentclass[11pt]{amsart} 
\usepackage{amsmath, amsthm, amscd} 
\usepackage{color,graphics,epic,psfig,epsfig} 
\input amssym.def
\input amssym.tex 
\newtheorem{thm}{Theorem}[section] 
\newtheorem*{thm*}{Theorem} 
\newtheorem{cor}[thm]{Corollary} 
\newtheorem*{cor*}{Corollary} 
\newtheorem{lem}[thm]{Lemma} 
\newtheorem{prop}[thm]{Proposition}

\newtheorem*{con*}{Conjecture} 
\newtheorem*{que*}{Question}
 
\theoremstyle{definition} 
\newtheorem{defn}[thm]{Definition}
 
\theoremstyle{remark} 
\newtheorem{rem}[thm]{Remark}

\newcommand{\G}{\mathcal{G}}

\newcommand{\bbC}{\mathbb{C}}
\newcommand{\bbZ}{\mathbb{Z}}
\newcommand{\bbR}{\mathbb{R}}
\newcommand{\bbN}{\mathbb{N}}

\begin{document} 

\title{Roundness properties of groups}
\author[Jean-Francois Lafont]{Jean-Fran\c{c}ois Lafont}
\address{Department of Mathematical Science,
Binghamton University,
Binghamton, NY 13902}
\email{jlafont@math.binghamton.edu}
\author[Stratos Prassidis]{Stratos Prassidis$^*$} 
\address{Department of Mathematics {\&} Statistics,
Canisius College,
Buffalo, NY 14208}
\email{stratos@canisius.edu}
\protect\thanks{$^*$ Partially Supported by a Canisius College Summer Research Grant}

\begin{abstract}
Roundness of metric spaces was introduced by Per Enflo as a tool to study uniform structures of
linear topological spaces. The present paper investigates geometric and topological properties
detected by the roundness of general metric spaces. In particular, we show that geodesic spaces 
of roundness $2$ are contractible, and that a compact Riemannian manifold with roundness $>1$ 
must be simply connected. We then focus our investigation on Cayley graphs of 
finitely generated groups.  One of our main results is that every Cayley graph of a free 
abelian group on $\geq 2$
generators has roundness $=1$.  We show that if a group has no Cayley graph of roundness $=1$,
then it must be a torsion group with every element of order $2,3,5$, or $7$.   
\end{abstract} 
\maketitle 

\section{Introduction}

In a series of papers Per Enflo (\cite{en1}, \cite{en4}, \cite{en3}) used the idea of
metric roundness to investigate the uniform structure of Banach spaces. Later the same idea
was used in \cite{pw} to compare uniform structures between normed and quasi-normed linear
topological spaces. An extension of this property ({\it generalized roundness}) was used
by Enflo in the solution of Smirnov's problem (\cite{en2}). Also, if a metric space
has non-trivial generalized roundness, then some positive power of the distance function 
is a {\it negative kernel} on the space (\cite{ltw}). Negative kernels on Cayley graphs of
discrete groups were used for proving the coarse Baum--Connes Conjecture (and thus the
Novikov Conjecture) for these groups (\cite{ghw}, \cite{gk}). 

We investigate the roundness and generalized roundness properties of general metric spaces.
The triangle inequality implies that any metric space has roundness at least $1$. Using the
results of Enflo, 
essentially, a metric space $X$ has roundness $p$, $1 < p \le 2$ if $p$ is the supremum of all $q$
so that quadrilaterals in $X$ are {\it thiner} 
than the ones in an $L_q$-space. With this in mind, our first result is not surprising.

\begin{thm*} Every CAT(0)-space has roundness $2$.
\end{thm*}

On the other hand, it should be noticed that there are CAT(0)-spaces whose generalized roundness
is equal to $0$ (\cite{farhar}).  All spaces with approximate midpoints have roundness 
between $1$ and $2$. It is reasonable to try to understand the extremal cases.

\begin{thm*}
Proper geodesic spaces that have roundness $2$ are
contractible. 
\end{thm*}

We also point out that, in section $1.19_+$ of \cite{gromov}, Gromov raises the
question of determining what types of spaces one can obtain when one
imposes a restriction on the distances achieved between all $r$-tuples of
points.  The previous theorem can be viewed as a partial answer to this
question in the context where the restrictions on the distances between
all $4$ -tuples of points are given by the roundness $= 2$ condition.

Combining these two results, we recover the well-known result that any proper CAT(0)-space is
contractible.

On the other hand,
it is interesting to notice that the roundness properties of 
metric spaces with non-trivial closed geodesics are very poor. Mild assumptions
on such a space imply that its roundness is $1$. In particular, we have:

\begin{thm*}
A non-simply connected, compact, Riemannian manifold has roundness $1$.
\end{thm*}

This is in fact a special case of a more general theorem applying to geodesic metric spaces
with non-trivial fundamental group, and satisfying an additional hypothesis on existence of
convex neighborhoods around every point.  This more general result suggests that, as far as 
roundness is concerned, the most interesting spaces to look at
are simply connected geodesic spaces or, at the other extreme, totally disconnected spaces.

Our main explicit calculations are on discrete metric spaces determined by graphs, in
particular Cayley graphs of finitely generated groups. Roundness is not a quasi-isometric
invariant and thus, in general, the roundness of a Cayley graph of a group depends on the
choice of generating set. 
%But we
%can determine some general properties of roundness of a group in special cases.
%\begin{thm*}
%\null
%\noindent
%\begin{enumerate}
%\item The roundness of any Cayley graph of an infinite group is between $1$ and $2$.
%\item A group has a presentation with Cayley graph of roundness $\infty$ if and only if it
%is finite.
%\item The finitely generated free group has two different Cayley graphs one of which
%has roundness $1$ and the other has roundness $2$.
%\end{enumerate}
%\end{thm*}
A more relevant algebraic invariant seems to be the {\it roundness spectrum} of a group, 
which is the
collection of the roundness of all the Cayley graphs of the group. One of our main results is:

\begin{thm*}
The roundness spectrum of a finitely generated free abelian group on more than one generator 
is $\{1\}$.
\end{thm*}

In general, the roundness spectrum has the following property:

\begin{thm*}
If the roundness spectrum of $G$ does not contain $1$ then $G$ is a purely torsion group
in which every element has order $2$, $3$, $5$ or $7$.
\end{thm*}

In \cite{ltw} it was shown that in spaces with generalized roundness $p > 0$ 
the $p$-th power of the distance function is a negative kernel. Using this result we show the
following:

\begin{thm*}
Let $G$ be a group having a presentation whose Cayley graph has positive generalized roundness.
Then $G$ satisfies the coarse 
Baum--Connes Conjecture and thus the strong Novikov Conjecture.
\end{thm*}

In particular, if $G$ be a group having a presentation whose Cayley graph
isometrically 
embeds into an  $L_p$-space with $1 \le p \le 2$, then $G$ satisfies the coarse 
Baum--Connes Conjecture and thus the strong Novikov Conjecture
(see also \cite{now}, Corollary 4.3).

On the other hand negative kernels are closely related to the Kazhdan property. Using this
we can show that:

\begin{thm*}
Every Cayley graph of a finitely generated infinite Kazhdan group has generalized roundness $0$.
\end{thm*}

This result follows from combining the fact that infinite Kazhdan groups do not admit 
negative kernels (\cite{dHV}, \cite{de}) 
and the equivalence between non-trivial generalized roundness and 
negative kernels (\cite{ltw}). It should be noted that generalized roundness is an easier 
condition to be checked than the existence of negative kernels because generalized roundness
is a property of finite subspaces of the space.

We would like to thank Tom Farrell, Ralf Spatzier and Tony Weston for 
their helpful suggestions during the
course of this work.

\section{Preliminaries}

\begin{defn}\label{roundness}
	Let $(X, d)$ be a metric space, $p \in [1, {\infty}]$.
\begin{enumerate}
\item The {\it roundness} of $(X, d)$ is $p$ if $p$ is the supremum of all  $q$ such that:
for any four points $x_{00}$, $x_{10}$, $x_{01}$, $x_{11}$ in $X$,
$$d(x_{00}, x_{11})^q + d(x_{01}, x_{10})^q\le 
d(x_{00}, x_{01})^q + d(x_{00}, x_{10})^q + d(x_{11}, x_{01})^q + d(x_{11}, x_{10})^q.$$
\item The {\it generalized roundness} of $(X, d)$ is the supremum of all $q$ such that:
for every $n \ge 2$ and any collection of $2n$-points $\{a_i\}_{i=1}^n$, $\{b_i\}_{i=1}^n$, we
have that:
$$\sum_{1\le i < j \le n}(d(a_i, a_j)^q + d(b_i, b_j)^q) \le
\sum_{1\le i,j\le n}d(a_i, b_j)^q.$$
\end{enumerate}
\end{defn} 

\begin{rem}
\null

\noindent
\begin{enumerate}
\item 
Definition \ref{roundness}, Part (1), can be rephrased in terms of $2$-cubes. Recall that the unit
cube in ${\bbR}^{n}$ ($n \in {\bbN}$) is the set of $n$-vectors $\{0,1\}^{n}$.
An $n$-\textit{cube} $N$ in an arbitrary metric space $(X,d)$ is a collection of $2^{n}$
(not necessarily distinct) points in $X$ where each point in the collection is indexed 
by a distinct $n$-vector $\varepsilon \in \{0,1\}^{n}$ from the unit cube. A \textit{diagonal}
in $N$ is a pair of vertices $(x_{\varepsilon},x_{\delta})$ such that $\varepsilon$ and $\delta$
differ in all coordinates. An \textit{edge} in $N$ is a pair of vertices $(x_{\varepsilon},
x_{\delta})$ such that $\varepsilon$ and $\delta$ differ in precisely one coordinate. The
set of diagonals in $N$ will be denoted $D(N)$ and the set of edges in $N$ will be
denoted $E(N)$. An $n$-cube $N$ has $2^{n-1}$ diagonals and $n2^{n-1}$ edges. If
$f=(x, y)$ is an edge or diagonal in $N$, we will let $l(f)$ denote the $d$-length of $f$ 
in $X$. In other words, $l(f) = d(x,y)$. The analytic condition in Definition \ref{roundness}, 
Part (1), 
is a statement about $2$-cubes $N$ in $X$:
\[
\sum\limits_{d \in D(N)} l(d)^{q} \leq
\sum\limits_{e \in E(N)} l(e)^{q}.
\]
\item
Enflo \cite{en1} showed that roundness has the following exceptionally nice inductive property:

\noindent
In an $n$-cube $N$ in a metric space $(X,d)$ with roundness $p$ we have:
\[
\sum\limits_{d \in D(N)} l(d)^{p} \leq \sum\limits_{e \in E(N)} l(e)^{p}.
\]
In particular, if $d_{\rm{min}}$ denotes a diagonal of minimal $d$-length in $N$ and 
$e_{\rm{max}}$ denotes an edge of maximal $d$-length in $N$, then
$l(d_{\rm{min}}) \leq n^{\frac{1}{p}} \cdot l(e_{\rm{max}})$.
\item The triangle inequality implies that any metric space has roundness $\ge 1$. If the
space has approximate midpoints, then its roundness is $\le 2$.
\item The collection of $2n$ points in the second part of the definition is usually called
an $n$-double simplex.
\end{enumerate}
\end{rem}

\section{Geometric aspects of roundness}

Roundness and curvature bounded from above are two metric properties. In this section, we
examine their connections.
Our first observation is that in spaces with complicated topology, roundness cannot be large.
We consider one of the simplest non-simply connected space first.

\begin{lem}\label{lem-circle}
The roundness of the circle is $1$.
\end{lem}

\begin{proof}
Let $x_{00}$, $x_{01}$, $x_{10}$, $x_{11}$ be four points on $S^1$ so that:
$$d(x_{00}, x_{01}) + d(x_{01}, x_{11}) = d(x_{00}, x_{11}), \quad
d(x_{01}, x_{11}) + d(x_{11}, x_{10}) = d(x_{01}, x_{10}).$$
Then, for $p > 1$,
$$\begin{array}{lll}
d(x_{00}, x_{11})^p + d(x_{01}, x_{10})^p  & = &
(d(x_{00}, x_{01}) + d(x_{01}, x_{11}))^p + (d(x_{01}, x_{11}) + d(x_{11}, x_{10}))^p \\
& > & d(x_{00}, x_{01})^p + d(x_{01}, x_{11})^p + d(x_{01}, x_{11})^p + d(x_{11}, x_{10})^p 
\end{array}$$
Thus the roundness can not be larger than $1$.
\end{proof}

\begin{rem}\label{rem-circle}
The generalized roundness of the circle is $1$: Lemma \ref{lem-circle} implies that
the generalized roundness is less than or equal to $1$. But in \cite{dela}, Theorem 6.4.5,
it is shown that $S^1$ isometrically embeds into an ${\ell}_1$-space, which has generalized
roundness $1$ (\cite{ltw}).
\end{rem}

\begin{prop}\label{prop-one}
Let $(X, d)$ be a geodesic metric space that admits a globally minimizing closed geodesic. Then
its roundness is $1$.
\end{prop}

\begin{proof}
A globally minimizing closed geodesic $\gamma$ is an isometric embedding of a circle of length
${\ell}({\gamma})$. The Proposition follows from Lemma \ref{lem-circle}.
\end{proof}

\begin{lem}\label{lem-inf}
Let $(X, d)$ be a geodesic space.
Suppose there is a closed curve $\gamma$ such that 
$${\ell}({\gamma}) = \inf\{{\ell}({\gamma}):\; \gamma
\;\text{homotopically non-trivial rectifiable curve}\} > 0.$$
Then the roundness of $X$ is equal to $1$.
\end{lem}

\begin{proof}
Assume that the roundness of $X$ is greater than $1$.
Proposition \ref{prop-one} implies that such $\gamma$ can not be 
globally length minimizing. Hence, if
$${\gamma}: [0, L] \to X$$
is a unit length parametrization we have, after reparametrizing if necessary, that there is
$s \in [0, L/2]$ such that $d({\gamma}(0), {\gamma}(s)) < s$. As $X$ is a geodesic space, there
is a curve $\eta$ from ${\gamma}(0)$ to ${\gamma}(s)$ whose length is equal to
$d({\gamma}(0), {\gamma}(s))$. Let ${\gamma}_1$ be $\gamma$ restricted to $[0, s]$,
${\gamma}_2$ be $\gamma$ restricted to $[s, L]$. Form two new loops:
$${\eta}_1 = {\eta}^{-1}*{\gamma}_1, \; {\eta}_2 = {\gamma}_2*{\eta}.$$
Note that ${\eta}_2*{\eta}_1 \simeq {\gamma}$. Since $\gamma$ represents a non-trivial element
in ${\pi}_1(X)$, one of the loops ${\eta}_1$, ${\eta}_2$ must likewise be non-trivial.
We now compute the lengths of ${\eta}_1$, ${\eta}_2$:
$$\begin{array}{lll}
{\ell}({\eta}_1) & = & {\ell}({\gamma}_1) + {\ell}({\eta}) = s + {\ell}({\eta}) < s + s = 2s \le
L \\
{\ell}({\eta}_2) & = & {\ell}({\gamma}_2) + {\ell}({\eta}) = (L - s) + {\ell}({\eta}) 
< (L - s) + s = L
\end{array}$$
So in both cases, we find a homotopically non-trivial loop with length shorter than the assumed
minimum $L$, contradiction.
\end{proof}

The above lemma can be applied to a certain natural class of metric spaces.

\begin{defn}
A metric space $(X, d)$ is called {\it good} provided that, for each $p \in X$, there is
a neighborhood $N_p$ of $p$ with:
\begin{enumerate}
\item $N_p$ is simply connected.
\item $N_p$ is geodesically convex i.e., for each $y, z \in N_p$ and for each
geodesic $\gamma$ joining $y$ to $z$ with ${\ell}({\gamma}) = d(y, z)$, the trace of
$\gamma$ is contained in $N_p$.
\end{enumerate}
\end{defn}

\begin{rem}\label{rem-riemannian}
If $(X, d)$ is a Riemannian manifold then $(X, d)$ is good; this follows from the existence
of normal neighborhoods.  More generally, any Finsler manifold is good (this is due to 
J.H.C. Whitehead \cite{wh}; the authors thank Z. Shen for informing us of this result).
\end{rem}

\begin{prop}\label{prop-good}
Let $(X, d)$ be a good, compact, geodesic space 
with non-trivial fundamental group. Then there is a loop $\gamma$ such that:
\begin{enumerate}
\item $\gamma$ is not freely homotopic to a constant loop.
\item For each loop ${\gamma}'$ not freely homotopic to a constant loop, ${\ell}({\gamma}') \ge
{\ell}({\gamma})$
\end{enumerate}
\end{prop}

\begin{proof}
Let $L = {\inf}\{{\ell}({\eta})| \; {\eta} \;\text{not freely homotopic to a constant loop}\}$ 
and let $\{{\gamma}_i\}_{i\in\bbN}$ be a sequence of loops, each of which is not freely homotopic
to a constant loop such that ${\ell}({\gamma}_i) \rightarrow L$. We first observe that, without
loss of generality, we can assume that ${\gamma}_i$ is piecewise geodesic. Indeed, for a given
${\gamma}_i$, we can cover the trace of ${\gamma}_i$ with a finite collection of simply connected,
geodesically convex neighborhoods $N_j$, $j = 1, \dots , k$, since $(X, d)$ is a good
geodesic space. Pick $t_j$, $j = 1, \dots , k$,
in $S^1$ such that ${\gamma}_j([t_j, t_{j+1}]) \subset N_j$ and replace ${\gamma}_j|[t_j, t_{j+1}]$
by a geodesic lying in $N_j$ and joining ${\gamma}_j(t_j)$ to ${\gamma}_j(t_{j+1})$. Since $N_j$
is simply connected, the new loop is freely homotopic to the original ${\gamma}_j$, is piecewise
geodesic, and it has length less than or equal to ${\ell}({\gamma}_j)$. Hence this new sequence
of loops also has lengths tending to $L$. 

Now parametrize each of these loops with respect to arclength, scaled by ${\ell}({\gamma}_i)$, and
let $M = {\sup}\{{\ell}({\gamma}_i):\; i\in \bbN\}$. Note that $M < \infty$, and that for all $i$,
all $x, y\in S^1$, we have :
$$d({\gamma}_i(x), {\gamma}_i(y)) \le {\ell}({\gamma}_i)d_{S^1}(x, y) \le Md_{S^1}(x, y).$$
Hence the family of curves $\{{\gamma}_i\}$ is equicontinuous, and as $X$ is compact, a 
subsequence (also denoted $\{{\gamma}_i\}$) converges to a closed loop ${\gamma}_{\infty}$.

\vspace{12pt}\noindent
{\bf Claim 1}. ${\gamma}_{\infty}$ is freely homotopic to ${\gamma}_i$, for sufficiently
large $i$.

\vspace{12pt}\noindent
{\bf Proof}. The assumptions on $X$ allow us to cover the trace of ${\gamma}_{\infty}$ by a finite
sequence of simply connected, geodesically convex neighborhoods $N_j$, $j = 1, \dots , k$.
As before, choose $t_j$, $j = 1, \dots , k$,
in $S^1$ such that ${\gamma}_j([t_j, t_{j+1}]) \subset N_j$. Note that, since ${\gamma}_i \to 
{\gamma}_{\infty}$ uniformly, we can also have that, for $i$ sufficiently large, that
${\gamma}_i([t_j, t_{j+1}]) \subset N_j$. For each $1 \le j \le k$, pick a geodesic ${\eta}_j$
joining ${\gamma}(t_j)$ to ${\gamma}_{\infty}(t_j)$. Let
$${\gamma}^j_i = {\gamma}_i|[t_j, t_{j+1}], \quad  {\gamma}^j_{\infty} = 
{\gamma}_{\infty}|[t_j, t_{j+1}].$$
Consider the closed loops $({\gamma}_j^i)^{-1}*({\eta}_{j+1})^{-1}*{\gamma}_{\infty}^j*{\eta}_j$,
and observe that this closed loop lies entirely in $N_j$. Since $N_j$ is simply connected, this
loop is contractible. Concatenating the homotopies on the various pieces, we see that 
${\gamma}_{\infty}$ is freely 
homotopic to ${\gamma}_i$, for $i$ sufficiently large, proving the claim.

\vspace{12pt}
Claim 1 implies:
\begin{enumerate}
\item ${\gamma}_{\infty}$ is not freely homotopic to a constant loop.
\item From the definition of $L$, we derive that ${\ell}({\gamma}_{\infty}) \ge L$.
\end{enumerate}
The rest of this proof is fairly standard.
Replace ${\gamma}_{\infty}$ by a curve $\gamma$ which is piecewise geodesic,
with geodesics joining ${\gamma}_{\infty}(t_j)$ and ${\gamma}_{\infty}(t_{j+1})$. As before,
${\gamma}$ is freely homotopic to ${\gamma}_{\infty}$, hence ${\ell}({\gamma}) \ge L$.

\vspace{12pt}\noindent
{\bf Claim 2}. ${\ell}({\gamma}) = L$.

\vspace{12pt}\noindent
{\bf Proof}. Assume not. Then ${\ell}({\gamma}) > L$. Let
$${\varepsilon} = \frac{{\ell}({\gamma}) - L}{2k + 1} > 0.$$
Then there exists a positive integer $i$ such that:
\begin{enumerate}
\item ${\ell}({\gamma}_j) - L < {\varepsilon}$.
\item $d({\gamma}_j(t), {\gamma}_{\infty}(t)) < {\varepsilon}$.
\end{enumerate}
As before, set 
$${\gamma}^j_i = {\gamma}_i|[t_j, t_{j+1}], \quad  {\gamma}^j = {\gamma}|[t_j, t_{j+1}].$$
We have that:
$$\sum_{j=1}^k({\ell}({\gamma}_i^j) + 2{\varepsilon}) =
{\ell}({\gamma}_i) + 2k{\varepsilon} < L + (2k + 1){\varepsilon} = {\ell}({\gamma}) =
\sum_{j=1}^k{\ell}({\gamma}^j).$$ 
Hence there is $j$ such that ${\ell}({\gamma}_i^j) + 2{\varepsilon} < {\ell}({\gamma}^j)$. But this
contradicts the fact that each ${\gamma}^j$ is a geodesic. Hence ${\ell} \le L$.
That completes the proof of Claim 2 and the proposition.
\end{proof}

Combining Lemma \ref{lem-inf}, Remark \ref{rem-riemannian} 
and Proposition \ref{prop-good} we have:
 
\begin{cor}\label{cor-good}
Let $(X, d)$ be a good, compact, geodesic space 
with non-trivial fundamental group. Then the roundness of $X$ is equal to $1$. In
particular, a compact non-simply connected Riemannian manifold has roundness $1$.
\end{cor}

\begin{rem}
Corollary \ref{cor-good} implies that, from the roundness point of view, the most
interesting Riemannian manifolds are the simply connected ones.
\end{rem}

\begin{prop}\label{prop-cat-round}
Let $(X, d)$ be a CAT(0)-space. Then $(X, d)$ has roundness $2$.
\end{prop}

\begin{proof}
Since CAT(0)-spaces have approximate midpoints (\cite{mh}, Proposition 1.11), the roundness
of $(X, d)$ is $\le 2$. Now we will show that the roundness is at least $2$. So let
$x_{00}$, $x_{01}$, $x_{10}$, $x_{11}$ be four points in $X$. Proposition 1.11 in \cite{mh} 
implies that there is a subembedding of the four points in ${\bbR}^2$. More precisely, there
are points $\overline{x}_{00}$, $\overline{x}_{01}$, $\overline{x}_{10}$, $\overline{x}_{11}$ 
in ${\bbR}^2$ such
that:
$$d(x_{ij}, x_{k{\ell}}) = d(\overline{x}_{ij}, \overline{x}_{k{\ell}})$$
whenever $(i, j)$ and $(k, {\ell})$ are different in one coordinate, and
$$d(x_{ij}, x_{k{\ell}}) \le d(\overline{x}_{ij}, \overline{x}_{k{\ell}})$$
whenever they differ in both coordinates. Thus:
$$\begin{array}{lll}
d(x_{00}, x_{11})^2 + d(x_{01}, x_{10})^2 & \le &
d(\overline{x}_{00}, \overline{x}_{11})^2 + d(\overline{x}_{01}, \overline{x}_{10})^2 \\
& \le & d(\overline{x}_{00}, \overline{x}_{01})^2 + d(\overline{x}_{00}, \overline{x}_{10})^2 +
d(\overline{x}_{11}, \overline{x}_{01})^2 + d(\overline{x}_{11}, \overline{x}_{10})^2 \\
& = & d(x_{00}, x_{01})^2 + d(x_{00}, x_{10})^2 +
d(x_{11}, x_{01})^2 + d(x_{11}, x_{10})^2 
\end{array}$$
The second inequality holds because ${\bbR}^2$, with the standard metric, has roundness $2$.
\end{proof}

Roundness $2$ imposes geometric and metric restrictions on the space.

\begin{prop}\label{prop-geodesic}
Let $(X, d)$ be a geodesic metric space of roundness $2$. 
For any two points $A$ and $B$ in $X$, there is
a unique geodesic connecting them.
\end{prop}

\begin{proof}
Assume that there are two geodesics between $A$ and $B$. Let $M_i$, $i = 1, 2$, be the midpoints
on the corresponding geodesics. Apply the roundness $2$ inequality:
$$|M_1M_2|^2 + |AB|^2 \le |AM_1|^2 + |M_1B|^2 + |AM_2|^2 + |M_2B|^2 = |AB|^2.$$
where $|xy|$ denotes the distance between the points $x$ and $y$.  The inequality above
immediately forces
$M_1 = M_2$. Iterating this procedure we see that the two geodesics coincide
on a dense set of points, so that by continuity, they must coincide.
\end{proof}

\begin{prop}\label{prop-wcont}
Let $(X, d)$ be a proper geodesic space such that any pair of points in $X$ can be joined by
a unique geodesic segment. Then $X$ is contractible.
\end{prop}

\begin{proof}
Let $p\in X$ be the base point, and let $I$ denote the interval $[0,1]$. 
Define $F: X\times I \to X$ by letting $F(q, t)$
to be the time-one reparametrization of the geodesic segment joining $q$ to $p$. To show that
$F$ is continuous, let $(q, t)$ be a point in $X{\times}I$ and 
$\{(q_n, t_n)\}_{n\ge 1}$  a sequence of points that converges to $(q, t)$. If $F$ fails to be
continuous at $(q, t)$, then there exists a subsequence, also denoted $\{F(q_n, t_n)\}_{n\ge 1}$ 
with
$$d(F(q_n, t_n), F(q, t)) \ge {\varepsilon}, \;\text{for all} \; n, \; \text{for some}\;
\varepsilon > 0.$$
We also obtain that, 
since $F(q_n, t_n)$ lies on a geodesic joining $q_n$ to $p$:
$$d(p, F(p_n, t_n)) \le \sup_n\{d(p, q_n)\}.$$
Since $\{q_n\}$ converges to $q$, the supremum on the right is bounded hence the points 
$F(q_n, t_n)$ lie in some closed ball of radius $R$ at $p$.
The properness of the metric of $X$ ensures that 
there is a convergent subsequence of $\{F(q_n, t_n)\}_{n\ge 1}$.
After re-parametrizing we assume that
$$\lim_{n\to\infty}F(q_n, t_n) = z \not= F(q, t).$$
Set $S = \{q_n\; | \; n \in \bbN\}$ with the metric induced from $X$.

\vspace{18pt}\noindent
{\bf Claim}. Under the above hypotheses,
$$\begin{array}{lll}
d(p, z) & = & d(p, F(q, t)) \\
d(q, z) & = & d(q, F(q, t))
\end{array}$$

\vspace{12pt}\noindent
{\bf Proof}. The continuity of the distance function implies that the 
function
$$d(p, -): S \to {\bbR}$$
is continuous. Notice that $d(p, F(q, t)) = td(p, q)$. 
Thus the continuity of multiplication implies that
$${\phi} = d(p, F(-, -)): S{\times}I \to {\bbR}$$
is also continuous. Continuity of $\phi$ along with the
fact that $\{(q_n, t_n)\}_{n\ge 1}$ converges to $(q, t)$ implies that
$$\begin{array}{rlll}
\displaystyle{\lim_{n\to\infty}{\phi}(q_n, t_n)} & = & {\phi}(q, t) & \Rightarrow \\
\displaystyle{\lim_{n\to\infty}d(p, F(q_n, t_n))} & = & d(p, F(q, t))  & \Rightarrow \\
\displaystyle{d(p, \lim_{n\to\infty}F(q_n, t_n))} & = & d(p, F(q, t)) & \Rightarrow \\
d(p, z) & = & d(p, F(q, t)).
\end{array}$$
As before, the function 
$${\psi} = d(q, F(-, -)): S{\times}I \to {\bbR}$$
is continuous. Then 
$$\begin{array}{rlll}
\displaystyle{\lim_{n\to\infty}{\psi}(q_n, t_n)} & = & {\psi}(q, t) & \Rightarrow \\
\displaystyle{\lim_{n\to\infty}d(q, F(q_n, t_n))} & = & d(q, F(q, t))  & \Rightarrow \\
\displaystyle{d(f(q), \lim_{n\to\infty}F(q_n, t_n))} & = & d(q, F(q, t)) & \Rightarrow \\
d(q, z) & = & d(q, F(q, t)).
\end{array}$$
This proves the claim.

\vspace{18pt}
Using the Claim, one can find a path $\eta$ joining $p$ to $q$ by concatenating the unique
geodesic
from $p$ to $z$ and the unique geodesic from $z$ to $q$. The Claim shows that the length of 
$\eta$ is
$${\ell}({\eta}) = d(p, F(q, t)) + d(F(q, t), q) = d(p, q).$$
The last equality follows because $F(q, t)$ is a point on the geodesic joining $p$ to $q$. Since
${\ell}({\eta})$ is equal to the distance between its two end-points, $\eta$ is a geodesic. Since
$z$ belongs to the unique geodesic from $p$ to $q$, $z$ must lie on $\eta$, and the Claim forces
$z = F(q, t)$. This contradicts the fact that $d(F(q_n, t_n), F(q, t)) \ge {\varepsilon}>0$, for
all $n$.
\end{proof}

\section{Roundness Properties of Groups}

In this section we look at the geometric properties of Cayley graphs of finitely generated groups. 
The graphs will
be considered as discrete metric spaces equipped with the combinatorial distance.
In the remainder of this paper, we will consider finite, symmetric (i.e., $g \in \Sigma 
\Rightarrow g^{-1}\in\Sigma$) generating sets $\Sigma$ which do not contain the identity. 
Note that if the group $G$ does not contain any elements of order $2$, then
the generating sets of $G$ have even cardinality.  For a 4-tuple of points $w,x,y,z$ in a metric
space $X$, we use the notation $[w,x,y,z]$ to denote the 1-double simplex whose diagonals
are $\{w,y\}$ and $\{x,z\}$.  By the roundness of a 1-double simplex we will mean the supremum of
exponents for which the roundness inequality holds {\it for that specific 1-double simplex}.  This
of course provides an upper bound for the roundness of the space $X$.  We will similarly
use the term generalized roundness of a $n$-double simplex to refer to the supremum of 
exponents for which the generalized roundness inequality holds for that specific $n$-double simplex.

The following is well-known (\cite{ns}, Proposition 2). We outline the proof for completeness.

\begin{prop}\label{prop-tree}
Let $X$ be an $\bbR$-tree. Then the roundness of $X$ is $2$.
\end{prop}

\begin{proof}
Geodesics in ${\bbR}$-trees have midpoints. So the roundness of $X$ is $\le 2$. 
Now, any four points in
an ${\bbR}$-tree have a convex hull as in Figure 1.

\begin{figure}[htbp]
\input{tree.pstex_t}
\caption{}
\end{figure}
There are two cases to be considered. One is the quadrilateral $[A, B, C, D]$ and the other
is the quadrilateral 
$[A, C, B, D]$. Direct calculation shows that in both cases the inequality holds for $p = 2$. It is
also easy to see that the only time that equality holds is if the points $A$, $B$, $C$ and $D$ 
are colinear in that order and $d(A, B) = d(C, D)$. Then the quadrilateral $[A, B, C, D]$ has
roundness $2$.
\end{proof}

\begin{cor}\label{cor-free}
The Cayley graph of a non-trivial 
free group with the standard set of generators has roundness $2$. Also
the Cayley graph of the free product of finitely many copies of the cyclic group of order $2$ has 
roundness $2$.
\end{cor}

\begin{rem}\label{rem-tree}
The generalized roundness of a tree is $\ge 1$: in \cite{dela}, Example 19.1.4, it is shown that
finite trees can be isometrically embedded into the cube of a finite ${\ell}_1$-space. Notice
that any $n$-double simplex in the tree will be embedded isometrically into an 
${\ell}_1$-space. Thus it will have generalized roundness $\ge 1$. Since the roundness of the
tree is $2$, the generalized roundness is between $1$ and $2$. 
\end{rem}

%\begin{cor}
%The generalized roundness of an infinite, locally finite tree $T$ is $1$, and 
%the generalized roundness of an $\bbR$-tree is $2$.
%\end{cor}

\begin{rem}\label{rem-quasi}
Roundness is not an invariant of quasi-isometries of metric spaces: let $Cay(F_2, \{x, y\})$ 
be the Cayley
graph of the standard presentation of the free group $F_2$ on two generators $x$ and $y$.
Then $\G$ is a tree and thus it has roundness $2$.
We will give a different presentation of $F_2$:
$$F_2 = {\langle}x, y, z_1, z_2, z_3, z_4:\; 
z_1 = x^{-1}y, \; z_2 = xy, \; z_3 = xy^{-1}, \;z_4 = x^{-1}y^{-1}{\rangle}.$$
Then in the new Cayley graph $Cay(F_2, {\Sigma})$ there is a quadrilateral as in Figure 2.

\begin{figure}[htbp]
\input{freen.pstex_t}
\caption{}
\end{figure}

The lengths of the sides is $1$ and the diagonals have length $2$. That implies that the 
roundness of $Cay(F_2, {\Sigma})$ is equal to $1$. 
But $Cay(F_2, \{x, y\})$ and $Cay(F_2, {\Sigma})$ 
are quasi-isometric as they are Cayley graphs of the same group.
\end{rem}

We suggest another invariant that comes closer into being a quasi-isometry invariant, at least
for infinite groups.

\begin{defn}
Let $G$ be a finitely presented discrete group. The {\it roundness spectrum} of $G$ is defined as
$${\rho}(G) = \{{\rho}(Cay(G, {\Sigma})):\; \Sigma \;\text{a generating set for}\; G\}.$$
\end{defn}

\begin{rem}

\null
\noindent
\begin{enumerate}
\item In general, ${\rho}(G) \subseteq [1, {\infty}]$.
\item In Remark \ref{rem-quasi}, we have shown that ${\rho}(F_2) \supseteq \{1, 2\}$. 
If we use the presentation:
$$F_2 = {\langle}x, y, z:\; 
z = y^{-1}x{\rangle},$$
then the roundness of the $Cay(F_2, {\Sigma})$ is $\le {\ln}3/{\ln}2$; the authors suspect that 
the previous inequality is actually an equality.
\end{enumerate}
\end{rem}

\begin{prop}\label{prop-spectrum}
Let $G$ be an infinite, finitely generated group. Then ${\rho}(G) \subset [1, 2]$.
\end{prop}

\begin{proof}
Let $\Sigma$ be a finite presentation of $G$. 
Assume that $Cay(G, \Sigma )$ contains three points, $x$, $y$,
and $z$, such that:
$$d(x, y) = d(y, z) = 1, \; d(z, x) = 2.$$
Then ${\rho}(Cay(G, \Sigma )) \le 2$ because $y$ is the midpoint of $x$ and $z$. 
If there is no such
triple, then the triangle inequality implies that, for all triples $x$, $y$, $z$,
$$d(x, y) = d(y, z) = 1 \;\Longrightarrow \; d(z, x) = 1.$$
Therefore, if $g$ and $h$ are generators so is $gh$. That implies $\Sigma = G$, a contradiction,
since $\Sigma$ is finite and $G$ is infinite.
\end{proof}

\begin{rem}
The previous result is not true for finite groups. For a finite group $G$, let $G$ be the
set of generators. The Cayley graph of this presentation is a finite complete graph.
But the roundness of a complete finite graph is $\infty$. So the roundness spectrum of a finite
group always contains $\infty$.
\end{rem}

Actually, the example in the Remark \ref{rem-quasi} suggests a way of constructing Cayley graphs 
for almost any group whose roundness is $1$.

\begin{prop}\label{prop-round_one}
Let $G$ be a finitely generated group, containing two elements $x$ and $y$ with the property that:
\begin{enumerate}
\item $x$ and $y$ do not have order $2$.
\item $x \not= y^{{\pm}1}.$
\item $x^3 \not= y^{{\pm}1}$ and $y^3 \not= x^{{\pm}1}$.
\end{enumerate}
Then $1\in \rho (G)$.
\end{prop}

\begin{proof}
Let ${\Sigma}$ be a finite symmetric set of generators of $G$.
Include $x$ and $y$ in the set of generators. If $x^2$ or $y^2$ belong to ${\Sigma}$, then 
remove those generators.
Also, include the generators and relations:
$$z_1 = x^{-1}y, z_2 = xy, z_3 = xy^{-1}, z_4 = x^{-1}y^{-1}.$$
Let $\G$ be the Cayley graph in the new presentation.
The quadrilateral $[x, y, x^{-1}, y^{-1}]$ has all vertices distinct (by (1) and (2)), and 
the edges all have length $1$, since we added $z_i$, $i = 1, \dots 4$, as generators. The
diagonals $d_1 = [x, x^{-1}]$ and $d_2 = [y, y^{-1}]$ have length two. That is because $x^2$ and
$y^2$ do not belong to the generating set and Conditions (2) and (3) ensure that $x^{{\pm}2}$,
$y^{{\pm}2}$ are not equal to $z_i$, $i = 1, \dots 4$.
This 4-point
configuration implies that the roundness is $1$.
\end{proof}

\begin{cor}\label{cor-2}
Assume $G$ is a finitely generated group with $1\notin \rho (G)$. Then $G$ is a torsion group
with every element of order $2$, $3$, $5$ or $7$.
\end{cor}

\begin{proof}
Let $g\in G$ have order $n$ bigger than or equal to $7$. A simple counting argument shows that
in ${\langle}g{\rangle}$, there exists an element $g'$ such that:
$$g' \not\in\{g, g^{n - 1}, g^3, g^{n - 3}\}$$
and $(g')^3 \not= g^{{\pm}1}$. Then the pair $\{g, g'\}$ satisfies the conditions of Proposition
\ref{prop-round_one}.

Let $G$ contain an element $g$ of order $4$. Then include $g$ in the generating set of $G$. If 
$g^2$ or  $g^4$ are in the generating set then delete them from
the generating set. Then
the quadrilateral $[1, g, g^2, g^3]$ has roundness $1$. That is because from the construction the
edges have length $1$ and the diagonals have length $2$. 

If $G$ contains an element $g$ of order $6$, include $g$ in the generating set. Delete any
generator from the original set which is a power of $g$. Then the pair
$\{g, g^2\}$ satisfies the conditions of Proposition
\ref{prop-round_one}.
\end{proof}

\begin{rem}
An argument identical to that of Proposition \ref{prop-round_one} and Corollary \ref{cor-2} 
implies that if a graph $\G$ has minimal cycles of length different from $3$,
$5$ and $7$, then the roundness of $\G$ is $1$.
\end{rem} 

Let ${\bbZ}^2$ denote the free abelian group on two generators. We will consider ${\bbZ}^2$ as
the integral lattice in ${\bbR}^2$ and we will use coordinates to denote elements of ${\bbZ}^2$.
Let $\overrightarrow{i} = (1, 0)$, $\overrightarrow{j} = (0, 1)$ denote the standard basis of
${\bbZ}^2$.

\begin{thm}\label{thm-z2}
If $\Sigma$ is a generating set for ${\bbZ}^2$, then the Cayley graph ${\G}_{\Sigma}$ of 
$({\bbZ}^2, {\Sigma})$ has roundness $1$. In other words, $\rho ({\bbZ}^2)=\{1\}$.
\end{thm}

\begin{proof}
Before starting the proof we make two simple observations:
\begin{enumerate}
\item If $\Sigma$ is a finite, symmetric generating set for the group $G$, and ${\phi}\in
\text{Aut}(G)$, then there is a canonical isometry  between $Cay(G, {\Sigma})$ and
$Cay(G, {\phi}({\Sigma}))$; in fact, $\phi$ induces the isometry.
\item If $\Sigma$ is a finite symmetric generating set for ${\bbZ}^2$, and there exist $g$ and $h$
in $\Sigma$ ($g \not= {\pm}h$) with $g {\pm} h \notin \Sigma$, then $[0, g, g+h, h]$ is a 4-tuple
with roundness equal to $1$.
\end{enumerate}

We will consider cases depending on  $|{\Sigma}|$. 

\vspace{12pt}\noindent
\underline{Case 1}. $|{\Sigma}| = 4.$

\noindent
Then ${\Sigma} = \{{\pm}u, {\pm}v\}$ where $u, v \in {\bbZ}^2$ are linearly independent.
Observation (2) immediately applies, hence the roundness of ${\G}_{\Sigma}$ equals 1.

\vspace{12pt}\noindent
\underline{Case 2}. 
$|{\Sigma}| = 6.$

\noindent
If $\Sigma$ is not of the form $\{{\pm}u, {\pm}v, {\pm}(u + v)\}$, then Observation (2) applies
and we are done. Hence assume that $\Sigma$ is of the form above, and observe that
$\{{\pm}u, {\pm}v\} \subseteq \Sigma$ is already a generating set for ${\bbZ}^2$. But we know that 
$\text{Aut}({\bbZ}^2)$ acts transitively on pairs of generating elements. Hence from Observation
(1), it is sufficient to compute the roundness of ${\G}_{\Sigma}$ where
$${\Sigma} = \{{\pm}\overrightarrow{i}, {\pm}\overrightarrow{j}, {\pm}(\overrightarrow{i}
+ \overrightarrow{j})\}.$$
We will show
that roundness is $1$ by contradiction. Assume that the roundness is equal to $p > 1$. Consider
the quadrilateral with vertices $(0, 0)$, $(0, 1)$, $(n, 0)$, $(n, 1)$. Then the roundness
condition reads:
$$(n + 1)^p + n^p \le n^p + n^p + 1^p + 1^p \;\Longrightarrow\;
(n + 1)^p - n^p \le 2.$$
But, taking limits, and noting that $p > 1$, 
$$\lim_{n\to\infty}[(n + 1)^p - n^p] = \infty.$$
So there is $n \in \bbN$, such that 
$$(n + 1)^p - n^p > 2.$$
Contradiction.

\vspace{12pt}\noindent
\underline{Case 3}. $|{\Sigma}| = 2k \ge 8.$

\noindent
Theorem \ref{thm-generating} (proved in the Appendix) implies that $\Sigma$ contains two
elements $u$ and $v$ ($u \not= {\pm}v$) with $u {\pm} v \notin \Sigma$. Hence
Observation (2) applies and we are done.

This concludes the argument for Theorem \ref{thm-z2}.
\end{proof}

The proof of Theorem \ref{thm-z2} can be modified to work for any finitely generated free abelian
group.

\begin{thm}\label{thm-zn}
If $\Sigma$ is a generating set for ${\bbZ}^n$ ($n \ge 2$), then the Cayley graph
${\G}_{\Sigma}$ of $({\bbZ}^n, {\Sigma})$ has roundness $1$.  In other words,
$\rho (\bbZ ^n)=\{1\}$ whenever $n\geq 2$.
\end{thm}

\begin{proof}
In Theorem \ref{thm-z2}, the case $n = 2$ has been dealt with, so we assume that $n \ge 3$. As it
was already observed in the proof of Theorem \ref{thm-z2}, if there is a pair $u$, $v$ in $\Sigma$
with $u \not= {\pm}v$ and $u {\pm}v\notin \Sigma$, the quadrilateral $[0, u, u + v, v]$ has
roundness $1$ forcing the roundness of the Cayley graph to be $1$. Hence if we have a generating
set $\Sigma$ such that the roundness of ${\G}_{\Sigma}$ is not $1$, then $\Sigma$ has the property:
$$\text{for each}\; u, v\in \Sigma, u \not= {\pm}v, \;\text{either}\; u + v \in {\Sigma} \;
\text{or}\; u - v\in \Sigma \quad (*)$$
If $u$ and $v$ are two linearly independent elements in $\Sigma$, they span a subgroup of 
${\bbZ}^n$ that is isomorphic to ${\bbZ}^2$. Furthermore, the set ${\Sigma} = {\Sigma}
\cap{\langle}u, v{\rangle}$ also satisfies property (*). But the argument in Theorem
\ref{thm-generating} 
shows that any generating set of ${\bbZ}^2$ having property (*) has cardinality $6$ and it has 
(up to relabeling) the form
$${\Sigma} = \{{\pm}u, {\pm}v, {\pm}(u {\pm} v)\}.$$
Hence the original generating set $\Sigma$ has the stronger property:
$$\text{for all}\; u, v \in {\Sigma}, \{u, v\}\;\text{linearly independent}, \;
\text{either}\; u + v \in {\Sigma} \;\text{or}\; u - v\in {\Sigma}\; \text{{\bf but not both}}. 
\quad 
(**)$$
Since $n \ge 3$, $\Sigma$ contains at least three linearly independent elements $u$, $v$ and $w$.
Using Property (**) we see that, up to relabeling, there are two possible cases:
\begin{itemize}
\item[]{\bf Case 1}. $u + v\in {\Sigma}$, $u + w\in {\Sigma}$, $v - w\in {\Sigma}$.
\item[]{\bf Case 2}. $u + v\in {\Sigma}$, $u + w\in {\Sigma}$, $v + w\in {\Sigma}$.
\end{itemize}
We now discuss each case separately.

\vspace{12pt}\noindent
{\bf Case 1}. Since $u + v \in \Sigma$, $w\in \Sigma$, Property (**) implies that either
$u + v + w\in \Sigma$ or $u + v - w\in \Sigma$ but not both.

Let us assume that $u + v + w\in \Sigma$. Since $u + w \in \Sigma$, $v \in \Sigma$, Property (**)
implies that $u - v + w\notin\Sigma$. Since $u\in\Sigma$, $v - w\in\Sigma$ property (*) again 
forces $u + v - w\in\Sigma$. But now we have that $u + v\in\Sigma$, $w \in \Sigma$ and both
$(u + v){\pm}w\in\Sigma$, contradicting (**).

On the other hand, if $u + v - w\in \Sigma$, since $u + v\in \Sigma$, $w \in \Sigma$, (**)
implies that $u + v + w\notin \Sigma$. As $u + w\in\Sigma$, $v\in\Sigma$, (**) forces
$u - v + w\in \Sigma$. But now we have $u \in \Sigma$, $v - w\in\Sigma$, and both
$u{\pm}(v - w)\in\Sigma$, contradicting (**). Thus Case 1 cannot occur.

\vspace{12pt}\noindent
{\bf Case 2}. In this case, we claim that the assumption implies that $\Sigma$ must contain
$nu + nv + (n - 1) w$, $nu + (n - 1)v + nw$, $(n - 1)u + nv + nw$ for infinitely many $n\in \bbN$.
If this were the case, linear independence of $u$, $v$ and $w$ implies that all these elements
are distinct, contradicting the finiteness of $\Sigma$.

We show the Claim by recursion on $n$. The fact that this triple of vectors with $n = 1$ lie in
the generating set follows from the hypotheses for 
Case 2. Notice that if $\Sigma$ contains $nu + nv + (n - 1) w$, $nu + (n - 1)v + nw$, and
$(n - 1)u + nv + nw$, then it must contain the elements
$$2nu + (2n - 1)v + (2n - 1)w, (2n - 1)u + 2nv + (2n - 1)w, (2n - 1)u + (2n - 1)v + 2nw.$$
To see this observe that (**) along with the hypotheses for Case 2 implies that $u - v$,
$u - w$ and $v - w$ are not in $\Sigma$. The hypotheses along with (**) and
$$[nu + nv + (n - 1)w] - [nu + (n - 1)v + nw] = v - w\notin\Sigma$$
implies that
$$[nu + nv + (n - 1)w] + [nu + (n - 1)v + nw] = 2nu + (2n - 1)v + (2n - 1)w\in\Sigma$$
One applies the same reasoning to obtain the other two elements. So we obtain that indeed:
$$
\left.
\begin{array}{lll}
nu + nv + (n - 1)w & \in & \Sigma \\
nu + (n - 1)v + nw & \in & \Sigma \\
(n - 1)u + nv + nw & \in & \Sigma
\end{array}\right\}
\;\Longrightarrow\;
\left\{
\begin{array}{lll}
2nu + (2n - 1)v + (2n - 1)w & \in & \Sigma \\
(2n - 1)u + 2nv + (2n - 1)w & \in & \Sigma \\
(2n - 1)u + (2n - 1)v + 2nw & \in & \Sigma
\end{array}\right.$$
But now for this second set of elements of $\Sigma$, we see that the differences are
$$u - v \not\in\Sigma, \; u - w \not\in\Sigma, \; v - w \not\in\Sigma $$
hence their sums must be in $\Sigma$; so we have:
$$
\left.
\begin{array}{lll}
nu + nv + (n - 1)w & \in & \Sigma \\
nu + (n - 1)v + nw & \in & \Sigma \\
(n - 1)u + nv + nw & \in & \Sigma
\end{array}\right\}
\;\Longrightarrow\;
\left\{
\begin{array}{lll}
(4n - 1)u + (4n - 1)v + (4n - 2)w & \in & \Sigma \\
(4n - 1)u + (4n - 2)v + (4n - 1)w & \in & \Sigma \\
(4n - 2)u + (4n - 1)v + (4n - 1)w & \in & \Sigma
\end{array}\right.$$
Finally, observe that for $n \in \bbN$, $4n - 1 > n$. We conclude that 
$$
\left.
\begin{array}{lll}
nu + nv + (n - 1)w & \in & \Sigma \\
nu + (n - 1)v + nw & \in & \Sigma \\
(n - 1)u + nv + nw & \in & \Sigma
\end{array}\right\} \quad
\text{for}\; n = 1, 3, 11, 43, 171, \dots$$
giving the desired contradiction in Case 2.

As we obtain a contradiction in all case, we conclude that there is no finite symmetric generating
set $\Sigma$ having property (*), and hence ${\G}_{\Sigma}$ has roundness $1$.
\end{proof}

\begin{cor}
Let $\Sigma$ be a finite generating set of ${\bbZ}^n$. Then the Cayley graph ${\G}_{\Sigma}$ has
generalized roundness $\le 1$.
\end{cor}

\begin{rem}
The attentive reader might wonder whether there is a simpler proof for Theorem \ref{thm-zn},
and indeed might be tempted to argue as follows.  Take a pair of linearly independant vectors
from the generating set for
${\bbZ}^n$, and consider the ${\bbZ}^2$ subgroup they generate.  From Theorem \ref{thm-z2},
this subgroup has generalized roundness =1, hence there are configurations in the subgroup
whose roundness is =1, which would force the roundness of ${\bbZ}^n$ to likewise be =1.  The
problem with this approach is that the distance on the ${\bbZ}^2$ subgroup induced by the
ambient ${\bbZ}^n$ might {\it not}, \`a priori, be isometric to a Cayley graph of ${\bbZ}^2$.  
In fact,
this approach can be tweaked to give an easy proof in most cases.  As long as there is a 
pair of linearly independant vectors $u,v\in \Sigma$ with the property that $|\langle u,v\rangle
\cap \Sigma|\neq 6$, the argument outlined above can be modified to work.
\end{rem}

\section{Generalized Roundness and Baum--Connes Conjecture}

Generalized roundness is connected with the existence of negative kernels which are used
in proving certain forms of the Baum--Connes Conjecture.

\begin{defn}
Let $X$ be a set. A real valued function $h$ on $X{\times}X$ is called a negative
kernel provided that:
\begin{enumerate}
\item $h(x, x) = 0$, for all $x\in X$.
\item $h(x, y) = h(y, x)$, for all $x, y \in X$.
\item For all $n$-tuples $x_1$, $x_2$, \dots , $x_n$ in $X$ and $a_1$, $a_2$, \dots , $a_n$ in
$\bbR$ satisfying $\sum_{i=1}^na_i = 0$, we have that
$$\sum_{i,j = 1}^n a_ia_jh(x_i, x_j) \le 0.$$
\end{enumerate}
\end{defn}

In \cite{ltw}, it was shown that:

\begin{prop}\label{prop-genround}
In a metric space $X$, the $p$-th power of the distance function is a negative kernel
if and only if it has generalized roundness $\geq p$.
\end{prop}

An immediate application of the above result is to the generalized roundness of Kazhdan groups
(\cite{dHV}, \cite{de}).

\begin{prop}\label{prop-kazhdan}
Let $\Sigma$ be a finite generating set for an infinite Kazhdan group $G$ and ${\G}_{\Sigma}$ the
corresponding Cayley graph. Then the generalized roundness of ${\G}_{\Sigma}$ is $0$.
\end{prop}

\begin{proof}
Assume that the generalized roundness of ${\G}_{\Sigma}$ is $p > 0$. Then by Proposition
\ref{prop-genround} we have that 
$$d^p_{\Sigma}: G{\times}G \to {\bbR}$$
is a negative kernel. Define ${\Phi}_p: G \to {\bbR}$ by:
$${\Phi}_p(g) = d^p_{\Sigma}(g, e).$$
Then by the left invariance of the metric on ${\G}_{\Sigma}$, we get that $d_{\Sigma}^p(x, y) = 
{\Phi}_p(x^{-1}y)$. Furthermore, observe that if $z_j \in \bbC$, $j = 1, \dots , n$
satisfy ${\sum}_{j=1}^nz_j = 0$, then for any collection of $n$ elements $g_j$ of $G$, 
an easy computation yields  %(writing $z_j = a_j + ib_j$)%
that: 
$$\sum_{j,k=1}^n z_j\overline{z_k}d^p_{\Sigma}(g_j, g_k) \leq 0.$$
%$$\begin{array}{lll}
%\displaystyle{\sum_{j,k=1}^n z_j\overline{z_k}d^p_{\Sigma}(g_j, g_k)} & = &
%\displaystyle{\sum_{j,k=1}^n (a_j + ib_j)\overline{(a_k + ib_k)}d^p_{\Sigma}(g_j, g_k)} \\[2ex] 
%& = &
%\displaystyle{\sum_{j,k=1}^n ((a_ja_k + b_jb_k) + (b_ja_k - a_jb_k)i)
%d^p_{\Sigma}(g_j, g_k)} \\[2ex] 
%& = &
%\displaystyle{\sum_{j,k=1}^n (a_ja_k + b_jb_k)d^p_{\Sigma}(g_j, g_k)} \\[2ex] 
%& = &
%\displaystyle{\sum_{j,k=1}^n a_ja_kd^p_{\Sigma}(g_j, g_k) + 
%\sum_{j,k=1}^n b_jb_kd^p_{\Sigma}(g_j, g_k)} \\[2ex] 
%& \le & 0.
%\end{array}$$
%where the last inequality comes from the fact that 
%$$\sum_{j=1}^n a_j = \sum_{j=1}^n b_j = 0$$
%and the fact that $d_{\Sigma}^p$ is a negative kernel.
Since $G$ is Kazhdan, this implies that ${\Phi}_p$ is bounded (see Delorme \cite{de}). 
But $p > 0$ and $G$ is infinite, hence we obtain a contradiction.
\end{proof}

To apply the above to the coarse Baum--Connes conjecture we need the following definition:

\begin{defn}
Let $X$, $Y$ be a pair of metric spaces. A (not necessarily continuous) map 
$f:X \to Y$ is a coarse embedding if there are non-decreasing
proper function ${\rho}_{\pm}:
[0, {\infty}) \to [0, {\infty})$ such that:
$${\rho}_-(d_X(x, y)) \le d_Y(x,y) \le {\rho}_+(d_X(x, y)), \quad \text{for all}\; x, y \in X,$$
and with $\lim _{t\to \infty} \rho_-(t)=\infty$.
Of particular interest is the case where $Y$ is a Hilbert space, with distance induced by 
the norm.  A discrete metric space $X$ is said to have bounded geometry, provided that for 
every $r>0$, there exists a uniform upper bound $N(r)$ on the cardinality of the metric 
balls of radius $r$.
\end{defn}

Note that a composition of coarse embeddings is still a coarse embedding.  Furthermore, if 
$\Gamma$ is a finitely generated group, then the identity map provides a coarse embedding from
any Cayley graph of $\Gamma$ to any other Cayley graph of $\Gamma$.  Hence if one Cayley graph
coarsely embedds into Hilbert space, they all coarsely embedd into Hilbert space.  In this 
situation we will say that the group $\Gamma$ coarsely embedds into Hilbert space, and ignore
any reference to a Cayley graph.

Now Yu (\cite{yu}) has shown that discrete metric spaces with bounded geometry that are 
coarsely embeddable into a Hilbert space satisfy the coarse Baum--Connes conjecture.  
In particular, since Cayley graphs of finitely generated groups have bounded geometry, if a
finitely generated group coarsely embedds into Hilbert space, then the coarse Baum--Connes 
Conjecture holds for the space, and hence the strong Novikov conjecture holds for the group 
in question (see \cite{yu}).  Recall that the strong Novikov conjecture asserts the injectivity of 
the classical assembly map for topological K-theory, and implies (amongst other things) the original
Novikov conjecture: that the higher signatures are homotopy invariants.

\begin{thm}\label{thm-bc}
Let $\Gamma$ be a finitely generated group, and assume that $\Gamma$ coarsely embedds into 
a metric space $X$ with generalized roundness $p>0$.  Then $\Gamma$ coarsely embedds 
in Hilbert space.  In particular, $\Gamma$ must satisfy the coarse Baum--Connes
conjecture, and hence the strong Novikov conjecture.
\end{thm}

\begin{proof}
We start by observing that, since $X$ has generalized roundness $p>0$, the $p^{th}$ power
of the distance function is a negative kernel.  Next we recall that a classic result of Schoenberg
\cite{schoen} states that given a negative kernel $h$ on a set
$X$, there exists a map $f:X\rightarrow \mathcal{H}$ into a Hilbert space $\mathcal{H}$ with
the property that $h(x,y)=||f(x)-f(y)||^2$.  So in our setting, there exists a map $f:X\rightarrow 
\mathcal{H}$ with the property that $d_X^p(x,y)=||f(x)-f(y)||^2$ for all $x,y\in X$.  
In particular, the map $f$ is a coarse embedding, with $\rho_-(t)=\rho_+(t)=t^{p/2}$.
Since $\Gamma$ coarsely embedds into $X$ by hypothesis, the composition yields the desired
coarse embedding into $\mathcal{H}$.  
\end{proof}

Two special cases are worth pointing out.  Note that an isometric embedding is a coarse embedding, 
and a quasi-isometric embedding is also a coarse embedding.  Furthermore, if a group acts 
properly discontinuously, cocompactly, freely, and isometrically, on a space $X$,
then $\Gamma$ and $X$ are quasi-isometric.  This immediately yields:

\begin{cor}\label{easycor}
Let $\Gamma$ be a finitely generated group, $X$ a metric space with generalized roundness $>0$, 
and assume that either of the following holds:
\begin{enumerate}
\item a Cayley graph of $\Gamma$ isometrically embedds into $X$, or
\item $\Gamma$ acts properly discontinuously, cocompactly, with finite stabilizers, by isometries
on $X$.
\end{enumerate}
Then $\Gamma$ is coarsely embeddable into Hilbert space.  In particular, $\Gamma$ must
satisfy the coarse Baum--Connes conjecture, and hence the strong Novikov conjecture.
\end{cor}

Note that a special case of the above corollary is the situation where some Cayley graph of 
$\Gamma$ has generalized roundness $>0$.  To obtain some further examples, we note that 
in \cite{ltw}, it was proved that the Banach spaces $L_p(\mu)$ (with $1\leq p\leq 2$) have
generalized roundness $\geq p$.  Hence we have:

\begin{cor}\label{thm-bc1}
Assume that the Cayley graph of a group $\Gamma$ admits an isometric embedding into an 
$L_p({\mu})$ space with $1 \le p \le 2$. Then $\Gamma$ satisfies the coarse Baum--Connes 
conjecture and thus the strong Novikov conjecture.
\end{cor}

We point out that a somewhat more general version of Corollary \ref{thm-bc1} can be found
in the work of Nowak \cite{now}.  We also mention that in the book by Deza-Laurent (\cite{dela}  
Chapter 19), conditions are given for graphs to be embeddable into an ${\ell}_1$-space. 
A natural question to ask is whether a converse to Corollary \ref{easycor} can hold.  Our next 
result is a partial counterexample to the converse:

\begin{prop}\label{counter}
There exists a group $\Gamma$ which is coarsely embeddable into Hilbert space, but fails to
satisfy the hypotheses in Corollary \ref{easycor}. 
\end{prop}

\begin{proof}
In Proposition \ref{prop-kazhdan}, we showed that all Cayley graphs of finitely generated
Kazhdan groups have generalized roundness $=0$.  In particular, if $\Gamma$ is a uniform lattice
in $Sp(n,1)$ or $F_{4(-20)}$, then $\Gamma$ is Kazhdan (see \cite{dHV}), and hence every 
Cayley graph of $\Gamma$ has generalized roundness $=0$.  This implies that $\Gamma$ cannot
be isometrically embedded into any space $X$ with generalized roundness $>0$, and hence
fails to satisfy hypothesis (1) in Theorem \ref{thm-bc1}.  

Next note that if $\Gamma$ satisfies hypothesis (2) in Theorem \ref{thm-bc1}, then picking 
a point $x\in X$, one can define a new distance $d_\Gamma$ on $\Gamma$ by setting 
$d_\Gamma(g,h):=d_X(g\cdot x,h \cdot x)$.  Note that this distance is left-invariant under the 
natural $\Gamma$ action on itself.  Furthermore, with this distance, the map $\phi:\Gamma\to X$ 
given by $\phi(g)=g\cdot x$ is an isometric embedding, and hence $d_\Gamma$ must have generalized
roundness $>0$.  But now the argument given in Proposition \ref{prop-kazhdan} applies verbatim
and yields a contradiction.

Finally, we note that $\Gamma$ acts isometrically on a quaternionic hyperbolic space or on the 
Cayley hyperbolic
plane, hence $\Gamma$ is $\delta$-hyperbolic.  But Sela \cite{sela} has proved that 
$\delta$-hyperbolic groups uniformly embedd into Hilbert space, giving the desired result.
\end{proof}

Let us point out that a consequence of work of Faraut-Harzallah \cite{farhar} implies that 
the generalized roundness of quaternionic hyperbolic spaces and of the Cayley hyperbolic plane 
is $=0$.  Note however that this does not, \`a priori, imply our Proposition \ref{prop-kazhdan}
for uniform lattices in $Sp(n,1)$ or $F_{4(-20)}$.  Indeed, the difficulty again lies in 
that generalized roundness is not well behaved with respect to coarse embeddings.  

We conclude this section by pointing out that Gromov \cite{gromov1} has
established the existence of finitely generated groups whose Cayley graph
{\it cannot} be uniformly embedded into Hilbert space.  An immediate consequence
of Corollary \ref{thm-bc1} is the following:

\begin{cor}
The groups constructed by Gromov in \cite{gromov1} cannot:
\begin{enumerate}
\item have a Cayley graph that isometrically embedds into a space of generalized roundness $>0$, 
\item act properly discontinuously, cocompactly, with finite stabilizers, by isometries
on a space with generalized roundness $>0$.
\end{enumerate}
\end{cor}

\section{Open Problems}

The calculations presented in this paper suggest a few of questions related to roundness and 
generalized roundness.

\begin{que*}
Is every CAT(0) space coarsely equivalent to a space with positive generalized roundness?
\end{que*}

Using Theorem \ref{thm-bc}, a positive answer to this question would imply
the coarse Baum--Connes Conjecture for groups acting properly discontinuously, freely and cocompactly
by isometries on $CAT(0)$-spaces. Note that while the Novikov Conjecture is known for these 
groups (\cite{cape},\cite{fala}), the coarse Baum--Connes is still open.

A well known result is that hyperbolic groups do not contain ${\bbZ}^2$ (\cite{mh}, Corollary 3.10). 
The following 
is an analogue for groups with non-trivial roundness spectrum.

\begin{que*}
Let $G$ be a finitely presented group such that its roundness spectrum is
strictly larger than $\{ 1\}$. Can $G$ contain a free abelian group of rank $2$?
\end{que*}

Concerning compact Riemannian manifolds, one can ask:

\begin{que*}
Does every compact Riemannian manifold contain a globally minimizing closed geodesic?
Do they always have roundness $=1$?
\end{que*}

We have answered both questions (see Proposition 3.7) for compact Riemannian 
manifolds with 
non-trivial fundamental group.  If the answer to the first question were affirmative in
general,
our Proposition 3.1 would immediately imply that all compact Riemannian manifolds have
roundness $=1$.  

In view of the fact that one of our main results is the computation of the roundness of Cayley
graphs of finitely generated free abelian groups, it is natural to ask:

\begin{que*}
What is the generalized roundness of a Cayley graph of ${\bbZ}^n$?
\end{que*}

It is clear from this paper that many of the difficulties in working with roundness and 
generalized roundness arise from the fact that these metric invariants are not coarse
invariants.  The authors believe that the development of coarse analogues of roundness 
and generalized roundness would be useful. The main hope would be that such a generalization 
would allow the results in Theorem \ref{thm-bc} to apply to a broader class of groups.

\section{Appendix}

We will show the combinatorial result used in the proof of Theorem \ref{thm-z2}. 
As before, let ${\bbZ}^2$ denote the free abelian group on two generators. 
Also, we embed ${\bbZ}^2$ as
the integral lattice in ${\bbR}^2$ and we will use coordinates to denote elements of ${\bbZ}^2$.
Let $\overrightarrow{i} = (1, 0)$, $\overrightarrow{j} = (0, 1)$ denote the standard basis of
${\bbZ}^2$. Let $\| - \|$ denote the usual norm on ${\bbR}^2$.
%$$\| \overrightarrow{a}, \overrightarrow{b} \| = \max\{|a_1 - b_1|, |a_2 - b_2|\}, 
%\;\text{where}\; \overrightarrow{a} = (a_1, a_2), \; \overrightarrow{b} = (b_1, b_2).$$

\begin{thm}\label{thm-generating}
Given a finite symmetric generating set $\Sigma$ with $|{\Sigma}| \ge 8$, 
then there is a pair
$g$ and $h$ in $\Sigma$, such that $g {\pm} h \notin \Sigma$.
\end{thm}

\begin{proof}
Let $\Sigma$ be a minimal generating set of cardinality bigger than or equal to $8$,
that satisfies property (*):
$$\text{for each}\; g, h\in \Sigma, g \not= {\pm}h, \;\text{either}\; g + h \in {\Sigma} \;
\text{or}\; g - h\in \Sigma \quad (*)$$
Then for any pair $\alpha$, $\beta$ of linearly independent elements of $\Sigma$ we have that
either
\begin{itemize}
\item[(1)] ${\langle}{\alpha}, {\beta}{\rangle} ={\bbZ}^2$ or
\item[(2)] $|{\Sigma}{\cap}{\langle}{\alpha}, {\beta}{\rangle}| = 6$.
\end{itemize} 
Indeed, if ${\langle}{\alpha}, {\beta}{\rangle}$ does not generate all of ${\bbZ}^2$, then
it generates a proper subgroup (isomorphic to ${\bbZ}^2$), and hence 
$|{\Sigma}{\cap}{\langle}{\alpha}, {\beta}{\rangle}|<|\Sigma|$.  But the subset 
$|{\Sigma}{\cap}{\langle}{\alpha}, {\beta}{\rangle}|$ is a generating set for the subgroup
${\langle}{\alpha}, {\beta}{\rangle}$ (which is abstractly a ${\bbZ}^2$), and inherits
the property (*).  By minimality of the cardinality of $\Sigma$, this implies that 
$|{\Sigma}{\cap}{\langle}{\alpha}, {\beta}{\rangle}| = 6$.

Notice that in case (2) above, we have that either 
${\alpha} + {\beta} \in {\Sigma}, \;\text{or}\;
{\alpha} - {\beta} \in {\Sigma},$
but {\bf not both}. We now break up the argument into cases.

\vskip 10pt

\par\noindent
(i) Assume that $\Sigma$ contains two elements that generate ${\bbZ}^2$. Then, after
applying an element of $SL(2,{\bbZ})$, we may assume that $\overrightarrow{i}$ and
$\overrightarrow{j}$ and $\overrightarrow{i} + \overrightarrow{j}$ are in $\Sigma$. 

\vskip 5pt

\par\noindent
(i-a) Assume that $\Sigma$ contains a vector
$\overrightarrow{v} = (v_1, v_2)$ such that $\min\{|v_1|, |v_2|\} > 2$.
We will show the proof when $\overrightarrow{v}$ is in the first quadrant. The other cases
follow similarly.
Notice that the pair $\{\overrightarrow{i}, \overrightarrow{v}\}$ is a linearly independent
subset but it does not generate ${\bbZ}^2$, since $|v-2| > 2$. The same true is for the pair
$\{\overrightarrow{j}, \overrightarrow{v}\}$, since $|v_1| > 2$. 
Thus the two pairs satisfy condition (2). Therefore
either $\overrightarrow{i} + \overrightarrow{v} \in {\Sigma}$ or
$\overrightarrow{i} - \overrightarrow{v} \in {\Sigma}$ but not both. 
Choose $\overrightarrow{v}$ to have
maximal norm among all elements of ${\Sigma}$ with both coordinates bigger than $2$. 
Since
$$\| \overrightarrow{v} + \overrightarrow{i}\| > \|\overrightarrow{v}\|, \;
\| \overrightarrow{v} + \overrightarrow{j}\| > \|\overrightarrow{v}\|$$ 
the maximality of $\| \overrightarrow{v}\|$ implies that
$$\overrightarrow{v} - \overrightarrow{i}\in {\Sigma}, \; 
\overrightarrow{v} - \overrightarrow{j}\in {\Sigma}.$$
Now consider the pair 
$\{\overrightarrow{v} - \overrightarrow{i}, \overrightarrow{v} - \overrightarrow{j}\}$. It is
a linearly independent subset and it does not generate ${\bbZ}^2$. Hence, either
$$(\overrightarrow{v} - \overrightarrow{i}) + (\overrightarrow{v} - \overrightarrow{j}) \in
{\Sigma}, \;\text{or}\;
(\overrightarrow{v} - \overrightarrow{i}) - (\overrightarrow{v} - \overrightarrow{j}) \in
{\Sigma},$$
but not both. Since
$$\|(\overrightarrow{v} - \overrightarrow{i}) + (\overrightarrow{v} - \overrightarrow{j})\|
> \|\overrightarrow{v}\|,$$
the maximality of $\|\overrightarrow{v}\|$ implies that
$$(\overrightarrow{v} - \overrightarrow{i}) - (\overrightarrow{v} - \overrightarrow{j}) =
\overrightarrow{j} - \overrightarrow{i} \in {\Sigma}.$$
Again $\{\overrightarrow{j} - \overrightarrow{i}, \overrightarrow{v}\}$ are linearly independent
and they do not generate ${\bbZ}^2$, so the sum or the difference, but not both are in ${\Sigma}$.
Assume, without loss of generality, that $\overrightarrow{v} + 
(\overrightarrow{j} - \overrightarrow{i}) \in {\Sigma}$. Then 
$$\overrightarrow{v} - 
(\overrightarrow{j} - \overrightarrow{i}) \notin {\Sigma}.\quad (*)$$
Since $\{\overrightarrow{v} - \overrightarrow{i}, \overrightarrow{j}\}$ are linearly independent,
do not generate ${\bbZ}^2$, and their sum is in ${\Sigma}$, their difference is not:
$$\overrightarrow{v} - 
\overrightarrow{i} - \overrightarrow{j} \notin {\Sigma}. \quad (**)$$
So $\{\overrightarrow{v} - \overrightarrow{j}, \overrightarrow{i}\}$ are linearly independent and
do not generate ${\bbZ}^2$ and
$$\begin{array}{lll}
(\overrightarrow{v} - \overrightarrow{j}) + \overrightarrow{i} & \notin {\Sigma} 
& \text{by $(*)$} \\
(\overrightarrow{v} - \overrightarrow{j}) - \overrightarrow{i} & \notin {\Sigma} 
& \text{by $(**)$} 
\end{array}$$
Contradiction.

\vskip 5pt

\noindent
(i-b) Assume that there is $\overrightarrow{v} = (v_1, v_2)$ in ${\Sigma}$ with $|v_1|$ maximal
and bigger than $2$. We assume $v_1 > 0$, the case $v_1 < 0$, follows from an identical argument.
Since $\overrightarrow{v} + \overrightarrow{j}$ or
$\overrightarrow{v} - \overrightarrow{j}$ belongs to ${\Sigma}$, we may assume $v_2 \not= 0$.
\par\noindent
{\bf Case 1}. Let $v_2 > 0$. Choose $\overrightarrow{v}$ so that $|v_1|$ is maximal and bigger than
$2$, and $v_2$ is positive and maximal. Then we have
$$\begin{array}{rlll}
\overrightarrow{v} - \overrightarrow{i} & \in & {\Sigma} &\text{(maximality of $v_1$, and
$v_1 > 0$)} \\
\overrightarrow{v} - \overrightarrow{j} & \in & {\Sigma} &\text{(maximality of $v_2$ amongst
$\overrightarrow{v}$ with $v_1$ maximal)} \\
\overrightarrow{v} - \overrightarrow{i} - \overrightarrow{j} & \in & {\Sigma} &
\text{(maximality of $v_1$ and
$\{\overrightarrow{v}, \overrightarrow{i} + \overrightarrow{j}\}\subset {\Sigma}$ 
with $v_1$ maximal)}
\end{array}$$
We also have:
$$\begin{array}{rlll}
\overrightarrow{v} - 2\overrightarrow{j} & \notin & {\Sigma} &\text{(
$\overrightarrow{v} - \overrightarrow{j}$ and $\overrightarrow{j}$ linearly
independent, do not generate ${\bbZ}^2$ and $\overrightarrow{v} \in {\Sigma}$)} \\
\overrightarrow{v} - \overrightarrow{i} - 2\overrightarrow{j} & \notin & {\Sigma} &
\text{($\overrightarrow{v} - \overrightarrow{i} - \overrightarrow{j}$ and
$\overrightarrow{j}$ are linearly independent,}\\ 
&&&\text{do not generate ${\bbZ}^2$ and $\overrightarrow{v} - 
\overrightarrow{i} \in {\Sigma}$)} \\
\end{array}$$
But since $\overrightarrow{v} - 
\overrightarrow{j} \in {\Sigma}$, $\overrightarrow{i} + 
\overrightarrow{j} \in {\Sigma}$ are linearly independent and do not generate ${\bbZ}^2$ we have
$$\begin{array}{rlll}
\text{either:} & (\overrightarrow{v} - 
\overrightarrow{j}) - (\overrightarrow{i} + 
\overrightarrow{j}) & = & \overrightarrow{v} - \overrightarrow{i} - 2\overrightarrow{j} \in 
{\Sigma} \\
\text{or:} & (\overrightarrow{v} - 
\overrightarrow{j}) + (\overrightarrow{i} + 
\overrightarrow{j}) & = & \overrightarrow{v} + \overrightarrow{i} \in 
{\Sigma}
\end{array}$$
But, as explained before, the first case could not occur. Thus 
$\overrightarrow{v} + \overrightarrow{i} \in {\Sigma}$, contradicting the maximality of $v_1$.

\noindent
{\bf Case 2}. We assume, as before, that $\overrightarrow{v} \in {\Sigma}$ with $v_1 > 0$ maximal
and bigger than $2$ and $v_2 < 0$. Then we have:
$$\begin{array}{rlll}
\overrightarrow{v} + \overrightarrow{j} & \in & {\Sigma} &\text{(minimality of $v_2$)} \\
\overrightarrow{v} - \overrightarrow{i} & \in & {\Sigma} &
\text{(maximality of $v_1$)} \\
\overrightarrow{v} - \overrightarrow{i} - \overrightarrow{j} & \in & {\Sigma} &\text{
(maximality of $v_1$ and}\; \{\overrightarrow{v}, \overrightarrow{i} + \overrightarrow{j}\}
\subseteq {\Sigma})
\end{array}$$
This forces $\overrightarrow{v} - \overrightarrow{j} \notin {\Sigma}$ (minimality of $v_2$)
and $\overrightarrow{v} - 2\overrightarrow{i} - \overrightarrow{j} \in {\Sigma}$ because 
$\{\overrightarrow{v} - \overrightarrow{i} - \overrightarrow{j}, \overrightarrow{i}\}$ are
linearly independent, do not generate ${\bbZ}^2$ and $\overrightarrow{v} 
- \overrightarrow{j} \notin {\Sigma}$. Hence both 
$$(\overrightarrow{v} - \overrightarrow{i}) {\pm} (\overrightarrow{i} + 
\overrightarrow{j}) \in {\Sigma}.$$
However, these are linearly independent and do not generate ${\bbZ}^2$, contradiction.

\vskip 5pt

\noindent
(i-c) The same argument shows that we can also exclude the case $|v_2| > 2$. Also, ${\Sigma}$ is
invariant under taking negatives. That is we need to exclude the following points:
$${\langle}2, {\pm}1{\rangle}, {\langle}2, {\pm}2{\rangle}, {\langle}2, 0{\rangle}, 
{\langle}1, 2{\rangle}, {\langle}1, -1{\rangle}, {\langle}2, 1{\rangle}, {\langle}0, 2{\rangle}, 
{\langle}-1, 2{\rangle}.$$ 
We just consider each case separately:
\par\noindent
(a) Assume that $\overrightarrow{v} = {\langle}2, 2{\rangle} \in {\Sigma}$. Then, 
${\langle}2, 1{\rangle}$ and ${\langle}1, 2{\rangle}$ are in ${\Sigma}$. Since the sum of the
last two vectors is not in ${\Sigma}$, thus
$${\langle}2, 1{\rangle} - {\langle}1, 2{\rangle} = {\langle}1, -1{\rangle} \in {\Sigma}.$$
But then $\overrightarrow{v} + {\langle}1, -1{\rangle}$ or
$\overrightarrow{v} + {\langle}1, -1{\rangle}$ must be in ${\Sigma}$. Contradiction because
one of the coordinates is greater than $2$.
\par\noindent
(b) Assume that $\overrightarrow{v} = {\langle}1, 2{\rangle} \in {\Sigma}$. Then,
${\langle}0, 2{\rangle}$ is in ${\Sigma}$. Since ${\langle}1, 1{\rangle}$ is in ${\Sigma}$,
${\langle}-1, 1{\rangle}$ is in ${\Sigma}$, which implies that 
${\langle}2, 1{\rangle}$ is in ${\Sigma}$. By (a), ${\langle}2, 1{\rangle} +
\overrightarrow{j}\not\in{\Sigma}$, we get ${\langle}2, 0{\rangle}$ is in ${\Sigma}$. Since
${\langle}1, 1{\rangle}$ and ${\langle}-1, 1{\rangle}$ are in ${\Sigma}$ and do not generate
${\bbZ}^2$ then only one of ${\langle}1, 1{\rangle}{\pm}{\langle}-1, 1{\rangle}$ can be
in ${\Sigma}$. That is a contradiction, because both ${\langle}2, 0{\rangle}$ and
${\langle}0, 2{\rangle}$ are in ${\Sigma}$.
\par\noindent
(c) In all the other cases, it is easy to see that  ${\langle}-1, 1{\rangle} \in {\Sigma}$.
By applying the matrix 
$$A = \left(
\begin{array}{rr}
-1 & 0 \\
0 & -1
\end{array}\right)$$
we can exclude ${\langle}-1, 2{\rangle}$ ${\langle}-2, 2{\rangle}$ (corresponding to cases
(a) and (b) above). That implies that ${\langle}0, 2{\rangle} \notin {\Sigma}$, because
$${\langle}0, 2{\rangle} {\pm} \overrightarrow{i} \notin{\Sigma}$$
from the previous cases. Also, ${\langle}2, 1{\rangle}$ can be excluded because
$${\langle}2, 1{\rangle} {\pm} \overrightarrow{j} \notin{\Sigma}.$$
All the other cases can be excluded similarly, except when ${\langle}-1, 1{\rangle}\in{\Sigma}$.
If ${\langle}-1, 1{\rangle}\in{\Sigma}$ then
$${\langle}-1, 1{\rangle} {\pm} {\langle}1, 1{\rangle} \notin{\Sigma},$$
from the previous cases.

\vskip 10pt

\noindent
(ii) Assume that no two elements of ${\Sigma}$ generate ${\bbZ}^2$. Let $\overrightarrow{x}$ and
$\overrightarrow{y}$ be two linearly independent elements of ${\Sigma}$. Then, after applying
an element of $SL(2, {\bbZ})$ we can assume that $m\overrightarrow{i}$ and
$n\overrightarrow{j}$ belong to ${\Sigma}$. If neither of the vectors 
$m\overrightarrow{i} \pm n\overrightarrow{j}$ are in $\Sigma$, then we are done.  
So let us assume that $m\overrightarrow{i} + n\overrightarrow{j}\in \Sigma$.

We now define three sets of points in $\mathbb Z^2$: 
\begin{itemize}
\item $\mathcal{L}_x$ consists of the integral points lying on the lines $y=0$, $y=\pm n$, $y=\pm 2n$,
\item $\mathcal{L}_y$ consists of the integral points lying on the lines $x=0$, $x=\pm m$, $x=\pm 2m$,
\item $\mathcal{L}_{xy}$ consists of the integral points lying on the lines $y=(n/m)x$, 
$y=(n/m)x \pm n$, $y= (n/m)x \pm 2n$.
\end{itemize}
Note that these subsets have the property that all of their pairwise intersections
lie in the subgroup of $\mathbb Z^2$ generated by the pair $(m\overrightarrow{i}, 
n\overrightarrow{j})$. This implies that the generating set $\Sigma$ must contain some 
vector $\overrightarrow{v}$ with the
property that: 
$$\overrightarrow{v}\notin (\mathcal{L}_x \cap \mathcal{L}_y)\cup 
(\mathcal{L}_x \cap \mathcal{L}_{xy}) \cup
(\mathcal{L}_y \cap \mathcal{L}_{xy}).$$ 
But now from basic set theory, we can conclude that this chosen vector also has the property
that:
$$\overrightarrow{v}\notin (\mathcal{L}_x \cup \mathcal{L}_y)\cap 
(\mathcal{L}_x \cup \mathcal{L}_{xy}) \cap
(\mathcal{L}_y \cup \mathcal{L}_{xy}).$$
Hence the vector $\overrightarrow{v}$ fails to lie in one of the pairwise intersection.  At the
cost of applying an automorphism of $\mathbb Z^2$, we may assume that we have a 
$\overrightarrow{v}\in \Sigma$ satisfying
$\overrightarrow{v}\notin \mathcal{L}_x \cup \mathcal{L}_y$.

But now observe that the argument in Case (i-a) works equally well in this setting.  Indeed, 
the fact that $\overrightarrow{v} \notin \mathcal{L}_x \cup \mathcal{L}_y$ implies that all
of the vectors $\overrightarrow{v} + (\epsilon_1 m\overrightarrow{i}) + (\epsilon_2
n\overrightarrow{j})$ are linearly independant from both $m\overrightarrow{i}$ and
$n\overrightarrow{j}$, where each $\epsilon_i\in \{0,\pm 1, \pm 2\}$.  In particular, carrying
out the argument in Case (i-a) but replacing each $\overrightarrow{i}, \overrightarrow{j}$
in that argument by $m\overrightarrow{i}, n\overrightarrow{j}$, we still have linear independance
at all the required steps.  Hence we again obtain a contradiction.  This concludes the proof
of Theorem 7.1.
\end{proof}

\end{document}